\documentclass{article}
\usepackage{amsthm}
\usepackage[T1]{fontenc}
\usepackage[frenchb, english]{babel}
\usepackage[hmargin=4cm,vmargin=4cm]{geometry}
\usepackage{amssymb}
\usepackage{amsmath}
\usepackage{pstricks}

\def\aro{{A}^{\!\!\!\raise5pt\hbox{$\scriptstyle \circ$}}}
\def\aroi{{A}^{\!\!\!\raise4pt\hbox{$\scriptscriptstyle \circ$}}}
\def\cro{\smash{{C}^{\!\!\!\raise5pt\hbox{$\scriptstyle \circ$}}}}
\def\croi{\smash{{C}^{\!\!\!\raise4pt\hbox{$\scriptscriptstyle \circ$}}}}

\begin{document}

\newtheorem{Th}[subsubsection]{Théorème}
\newtheorem{Pro}[subsubsection]{Proposition}
\newtheorem{De}[subsubsection]{Définition}
\newtheorem{Prt}[subsubsection]{Propriété}
\newtheorem{Prts}[subsubsection]{Propriétés}
\newtheorem{Le}[subsubsection]{Lemme}
\newtheorem{Hyp}[subsubsection]{Hypothèse}
\newtheorem{Cor}[subsubsection]{Corollaire}



\title{Cramér's theorem\\
for asymptotically decoupled fields}
\date{March 17, 2011}
\author{Raphaël Cerf and Pierre Petit\\ \\
Universit\'e Paris Sud}
\maketitle

\selectlanguage{english}

\begin{abstract}
We give a general setting for Cramér's large deviations theorem for the empirical means $\mathfrak{m}_{\Lambda(n)} \sigma$ of a field $\sigma$ of random vectors indexed by $\mathbb{Z}^d$, which contains Cramér's theorem for i.i.d. random vectors of \cite{Pet11a} and Sanov's theorem for asymptotically decoupled measures showed in \cite{Pfi01}. We add the notions of \emph{local control} and \emph{local convex tension} so that the general proofs in the i.i.d. case adapt to the asymptotically decoupled case. We obtain the standard weak large deviations principle for the sequence $(\mathfrak{m}_{\Lambda(n)} \sigma)_{n \geqslant 1}$ and the identification between the entropy and the opposite of the Fenchel-Legendre transform of the pressure, $s = -p^*$.
\end{abstract}

\selectlanguage{french}

\begin{abstract}
Nous établissons un cadre général pour le théorème de Cramér sur les grandes déviations des moyennes empiriques $\mathfrak{m}_{\Lambda(n)} \sigma$ d'un champ $\sigma$ de vecteurs aléatoires indexés par $\mathbb{Z}^d$, cadre qui contient le théorème de Cramér pour des vecteurs aléatoires i.i.d. de \cite{Pet11a} et le théorème de Sanov pour les mesures asymptotiquement découplées montré dans \cite{Pfi01}. Nous ajoutons les notions de \emph{contrôle local} et de \emph{convexe-tension locale}, de sorte que les preuves du cas i.i.d. s'adaptent au cas asymptotiquement découplé. On obtient le principe de grandes déviations faible classique pour la suite $(\mathfrak{m}_{\Lambda(n)} \sigma)_{n \geqslant 1}$ ainsi que l'identification entre l'entropie et l'opposée de la transformée de Fenchel-Legendre de la pression, $s = -p^*$.
\end{abstract}

\medskip

\textbf{Remarque :} La prochaine version de ce texte sera en langue anglaise.

\section{Introduction}

Le théorème de Cramér (dans sa version la plus générale) concerne les grandes déviations des moyennes empiriques $\overline{X}_n$ d'une suite $(X_n)_{n \geqslant 1}$ de variables aléatoires i.i.d. à valeurs dans un espace vectoriel localement convexe (cf. \cite{Pet11a}). Le présent texte étend la théorie de Cramér en relaxant l'hypothèse d'indépendance. Nous importons ici les idées du cas indépendant dans le cadre proposé par Pfister \cite{Pfi01}. Son approche des théorèmes de Sanov (au niveau 3, dans la terminologie de \cite{Ell85}) pour les mesures de Gibbs met en lumière la sous-additivité sous-jacente, donnant un cadre général englobant les résultats précédents (\cite{Geo93}, mais aussi \cite{DoV75} pour les chaînes de Markov). Voyant, comme \cite[chapitre 24]{Cer07}, le théorème de Sanov comme une conséquence du théorème de Cramér, nous généralisons donc la théorie de Cramér pour les variables indépendantes au cadre asymptotiquement découplé de \cite{Pfi01}.

\medskip

Le cadre contient le cas indépendant, notamment les versions de \cite{BaZ79}, de \cite{Cer07} et de \cite{Pet11a}. La question de la mesurabilité des moyennes empiriques est résolue en considérant des tribus initiales (qui sont égales aux tribus produits dans le cas indépendant). Ensuite, nous complétons l'hypothèse de \emph{découplage asymptotique} de Pfister \cite{Pfi01} (assurant un analogue de la propriété \textsf{(SA)} du lemme sous-additif \ref{lsa2}) par une hypothèse de \emph{contrôle local} (pour avoir l'analogue de la propriété \textsf{(C)} de \ref{lsa2}) et par la \emph{convexe-tension locale} (analogue de la convexe-tension dans le cas i.i.d.). Les deux premières hypothèses permettent donc d'établir un lemme sous-additif (ou plutôt pseudo-sous-additif) duquel découle le PGD faible, comme dans le cas indépendant. La suite de la démonstration est analogue à celle du cas indépendant. On notera toutefois qu'on n'a plus \textsf{(BS$_c$)} mais une forme affaiblie.

\section{Cadre}\label{scadreadi}

Après quelques préliminaires, nous introduisons deux notions qui permettrons d'établir un lemme sous-additif pour des champs non indépendants analogue au lemme clef sur lequel reposait la démonstration dans le cas i.i.d. : le découplage asymptotique inférieur et le contrôle local. Ces deux notions correspondent, dans le cas i.i.d., aux propriétés de sous-additivité \textsf{(SA)} et de contrôle \textsf{(C)}. D'autre part, l'invariance par translation sur le réseau $\mathbb{Z}^d$ remplacera l'identique distribution.

\subsection{Espace des configurations}

Soient $X$ un espace vectoriel réel, $\mathcal{F}$ une tribu sur $X$ stable par dilatation et translation et $d$ un entier strictement positif. On appelle \emph{boîte} un sous-ensemble cubique de $\mathbb{Z}^d$, autrement dit un translaté de
$$
\Lambda(n) := [0, n[^d \cap \mathbb{Z}^d
$$
pour un entier $n \geqslant 1$. On notera $\mathbb{B}$ l'ensemble des boîtes de $\mathbb{Z}^d$. Si $\smash{\sigma \in X^{\mathbb{Z}^d}}$ et $z \in \mathbb{Z}^d$, on notera $\sigma(z)$ la coordonnée suivant $z$ de $\sigma$. Pour tout $\Lambda \in \mathbb{B}$, on définit l'application
$$
\mathfrak{m}_\Lambda : \sigma \in X^{\mathbb{Z}^d} \mapsto \frac{1}{|\Lambda|} \sum_{z \in \Lambda} \sigma(z) \in X
$$
Pour tout $S$ sous-ensemble fini de $\mathbb{Z}^d$, on définit, sur $\smash{X^{\mathbb{Z}^d}}$, la tribu $\mathcal{F}_S$ initiale pour les applications
$$
\mathfrak{m}_\Lambda : X^{\mathbb{Z}^d} \rightarrow X
$$
où $\Lambda$ décrit l'ensemble des boîtes incluses dans $S$. La tribu $\mathcal{F}_S$ est donc la plus petite tribu sur $X^{\mathbb{Z}^d}$ qui rende ces applications mesurables. Notons que, si $S_1 \subset S_2$, alors $\mathcal{F}_{S_1} \subset \mathcal{F}_{S_2}$. D'autre part, si $S \in \mathbb{B}$ et $z \in \mathbb{Z}^d$, alors $\mathcal{F}_{S + z} = \{ A + z \, ; \, A \in \mathcal{F}_S \}$ où, pour tout $\smash{A \subset X^{\mathbb{Z}^d}}$, $A + z = \{ \sigma(\cdot - z) \, ; \, \sigma \in A \}$.

\medskip

\textbf{Remarque :} Pour tout $S$ sous-ensemble fini de $\mathbb{Z}^d$, on a toujours l'inclusion $\smash{\mathcal{F}^{\otimes S} \subset \mathcal{F}_S}$. Etant donné que $\mathcal{F}$ est stable par dilatation, l'égalité a lieu si et seulement si l'addition vectorielle
$$
(x, y) \in (X^2, \mathcal{F}^{\otimes S}) \mapsto x + y \in (X, \mathcal{F})
$$
est mesurable. C'est le cas si $\mathcal{F}$ est la tribu initiale pour une famille $(N_i, f_i)_{i \in I}$ où, pour tout $i \in I$, $N_i$ est un espace vectoriel normé séparable muni de sa tribu borélienne. Mais, étant donné que, pour tout $\Lambda \in \mathbb{B}$, nous allons nous donner une loi jointe sur $\mathcal{F}_\Lambda$, et non simplement, comme dans le cas i.i.d., une loi sur $\mathcal{F}$ qui permet de construire un $|\Lambda|$-uplet i.i.d., nous n'avons pas besoin de l'égalité $\smash{\mathcal{F}^{\otimes \Lambda} = \mathcal{F}_\Lambda}$.

\medskip

Enfin, pour tout $\Lambda \in \mathbb{B}$, on se donne une mesure de probabilité $\mu_\Lambda$ sur $\mathcal{F}_\Lambda$. On dira que $(\mu_\Lambda)_{\Lambda \in \mathbb{B}}$ est un \emph{système compatible de probabilités invariantes par translation} (\emph{s.c.i.t.}) sur $\smash{X^{\mathbb{Z}^d}}$ si :

\medskip

\textsf{(COMP)} pour toutes boîtes $\Lambda_1$ et $\Lambda_2$ de $\mathbb{Z}^d$ avec $\Lambda_1 \subset \Lambda_2$, on a :
$$
\forall A \in \mathcal{F}_{\Lambda_1} \qquad \mu_{\Lambda_2}(A) = \mu_{\Lambda_1}(A)
$$

\medskip

\textsf{(INV)} pour toute boîte $\Lambda$ de $\smash{\mathbb{Z}^d}$ et pour tout $z \in \mathbb{Z}^d$, on a :
$$
\forall A \in \mathcal{F}_\Lambda \qquad \mu_{\Lambda + z}(A + z) = \mu_\Lambda(A)
$$

\medskip

\textbf{Remarque :} Dans de nombreux cas étudiés (cas indépendant, mesures de Gibbs, chaînes de Markov...), il existe une mesure en volume infini décrivant les interactions, \emph{i.e.} une probabilité $\mathbb{P}$ sur $\smash{X^{\mathbb{Z}^d}}$ (c'est dans ce cadre que se place \cite{Pfi01}). Dans ce cas, le s.c.i.t. $(\mu_\Lambda)_{\Lambda \in \mathbb{B}}$ est simplement la famille des lois marginales de $\mathbb{P}$ sur les $\smash{X^\Lambda}$, pour $\Lambda \in \mathbb{B}$ (les deux conditions \textsf{(COMP)} et \textsf{(INV)} sont bien vérifiées). Cependant, il n'existe pas toujours de mesure en volume infini. Cependant, l'étude des grandes déviations ne nécessite que les marginales fini-dimensionnelles. Nous avons donc introduit ici la notion de s.c.i.t. pour avoir une hypothèse minimale. Toutefois, la construction habituelle de mesures en volume infini (les mesures de Gibbs, par exemple) repose sur un système de probabilités conditionnelles compatible au sens de Dobrushin, Lanford et Ruelle (DLR), et non au sens de Kolmogorov comme les systèmes que nous avons dénommés s.c.i.t. Ainsi, au lieu de construire une mesure en volume infini avec un système DLR, puis de définir le s.c.i.t. des marginales associées, il est sans doute possible d'obtenir un principe de grandes déviations directement sur le système DLR de départ. D'autant que, dans le cas indépendant, le système DLR et le s.c.i.t. coïncident. Ce sera l'objet d'un travail ultérieur.

\medskip

Etant donné un tel s.c.i.t. $(\mu_\Lambda)_{\Lambda \in \mathbb{B}}$, on définit, pour tout $\Lambda \in \mathbb{B}$ et pour tout $A \in \mathcal{F}_\Lambda$,
$$
\mathbb{P}(A) = \mu_\Lambda(A)
$$
L'hypothèse \textsf{(COMP)} assure que $\mathbb{P}(A)$ est bien défini pour tout $A$ dans l'union des $\mathcal{F}_\Lambda$ où $\Lambda \in \mathbb{B}$ (c'est une mesure additive sur le clan engendré par les $\mathcal{F}_\Lambda$, pour $\Lambda \in \mathbb{B}$). On définit également, pour $\Lambda \in \mathbb{B}$, $A \in \mathcal{F}_\Lambda$ et $\mathcal{G}$ une sous-tribu de $\mathcal{F}_\Lambda$,
$$
\mathbb{P}(A | \mathcal{G}) = \mu_\Lambda (A | \mathcal{G})
$$
et, pour $\Lambda \in \mathbb{B}$ et $\smash{Z : X^{\mathbb{Z}^d} \rightarrow [-\infty , +\infty]}$ une application $\mathcal{F}_\Lambda$-mesurable dont l'intégrale contre $\mu_\Lambda$ a un sens,
$$
\mathbb{E}(Z) = \int Z(\sigma) d\mu_\Lambda (\sigma)
$$
Si $D$ est une partie de $\smash{X^{\mathbb{Z}^d}}$, on dira que $(\mu_\Lambda)_{\Lambda \in \mathbb{B}}$ est \emph{porté par $D$} si, pour tout $\Lambda \in \mathbb{B}$, pour tout $(A, B) \in (\mathcal{F}_\Lambda)^2$,
$$
A \cap D = B \cap D \Rightarrow \mu_\Lambda(A) = \mu_\Lambda(B)
$$
C'est le cas si et seulement si, pour tout $\Lambda \in \mathbb{B}$, $\mu_\Lambda$ est la loi d'un champ $\sigma$ à valeurs dans $D$ (cf. cas i.i.d.). Comme dans le cas i.i.d., on fera souvent appel à cette description. Remarquons que $\sigma$ dépend de $\Lambda$, mais nous ne préciserons pas cette notation.

\subsection{Découplage asymptotique}

Soient $X$ un espace vectoriel réel, $\mathcal{F}$ une tribu sur $X$ stable par dilatation et translation, $d$ un entier strictement positif et $(\mu_\Lambda)_{\Lambda \in \mathbb{B}}$ un s.c.i.t. sur $\smash{X^{\mathbb{Z}^d}}$. Nous reprenons ici la définition de découplage de \cite{Pfi01}. Elle donne une condition de dépendance faible entre ce qui se passe dans une boîte et loin de cette boîte. Explicitons. On notera $\textup{dist}$ la distance associée à la norme $|\cdot|_\infty$ sur $\mathbb{Z}^d$. On dit que $(\mu_\Lambda)_{\Lambda \in \mathbb{B}}$ est \emph{asymptotiquement découplé inférieurement} (\emph{a.d.i.}) s'il existe deux applications $g$ et $c$ de $\mathbb{N}^*$ dans $[0, + \infty[$ telles que
$$
\frac{g(m)}{m} \to 0 \qquad \textrm{et} \qquad \frac{c(m)}{|\Lambda(m)|} \to 0
$$
et telles que : pour tout $m \geqslant 1$, pour tout sous-ensemble fini $S$ de $\mathbb{Z}^d$, pour tous $A \in \mathcal{F}_{\Lambda(m)}$ et $B \in \mathcal{F}_S$,
$$
\textup{dist} \big( S , \Lambda(m) \big) > g(m)
\Rightarrow \mathbb{P}(A \cap B) \geqslant e^{-c(m)} \mathbb{P}(A) \mathbb{P}(B)
$$
On dit que $(g, c)$ est un \emph{paramètre de découplage} de $(\mu_\Lambda)_{\Lambda \in \mathbb{B}}$. On verra par la suite que nous suffirait l'inégalité pour des $A$ et $B$ convexes mesurables particuliers. En effet, on aura seulement besoin de : pour tout $k \geqslant 1$, pour toutes boîtes $\Lambda_1$, ..., $\Lambda_k$ de $\mathbb{Z}^d$ de taille $m$ et telles que, pour tous $1 \leqslant i < j \leqslant k$, $\textup{dist}(\Lambda_i, \Lambda_j) > g(n)$, pour tous $C_1$, ... , $C_k$ translatés de convexes internes,
$$
\mathbb{P}\left( \bigcap_{i=1}^k \{ \mathfrak{m}_{\Lambda_i} \in C_i \} \right) \geqslant e^{-(k-1)c(n)} \prod_{i=1}^k \mathbb{P} (\mathfrak{m}_{\Lambda_i} \in C_i)
$$
Mais, en pratique, l'inégalité de découplage sera vérifiée pour tous les $A$ et $B$ mesurables.




\subsection{Contrôle local}

On reprend les notations de la section précédente. Si $C$ est un convexe de $X$, sa jauge (ou fonctionnelle de Minkowski) est l'application $M_C : X \rightarrow [0, +\infty]$ définie par
$$
\forall x \in X \qquad M_C(x) = \inf \{ t \geqslant 0 ; x \in tC \}
$$
On dit que $C$ est \emph{interne} si sa jauge $M_C$ est finie partout et si
$$
C = \{ x \in X \, ; \, M_C(x) < 1 \}
$$
L'ensemble des convexes internes est l'ensemble des convexes, voisinages ouverts de $0$ pour une certaine topologie localement convexe sur $X$. On note $\mathcal{C}_0(\mathcal{F})$ l'ensemble des convexes internes mesurables de $X$. On veut maintenant définir une condition de contrôle de chaque site conditionnellement au reste de la configuration. Soit $(\mu_\Lambda)_{\Lambda \in \mathbb{B}}$ un s.c.i.t. sur $\smash{X^{\mathbb{Z}^d}}$. On dit que $(\mu_\Lambda)_{\Lambda \in \mathbb{B}}$ est \emph{contrôlé localement} (\emph{c.l.}) s'il existe deux applications $t$ et $\alpha$ de $\mathcal{C}_0(\mathcal{F})$ dans $]0, +\infty[$ telles que : pour tout $C \in \mathcal{C}_0(\mathcal{F})$ et pour tout $\Lambda \in \mathbb{B}$,
$$
\mathbb{P} \big( \sigma(0) \in t(C) \cdot C | \mathcal{F}_{\Lambda \setminus \{ 0 \}} \big) \geqslant \alpha(C)
$$
On dit que $(t, \alpha)$ est un \emph{paramètre de contrôle} de $(\mu_\Lambda)_{\Lambda \in \mathbb{B}}$. Comme pour le découplage asymptotique inférieur, on pourrait se contenter de l'inégalité pour des conditionnements par des ensembles $A$ de la forme
$$
\bigcap_{i=1}^k \{ \mathfrak{m}_{\Lambda_i} \in C_i \}
$$
avec $\Lambda_1$, ..., $\Lambda_k$ boîtes disjointes de $\mathbb{Z}^d \setminus \{ 0 \}$ et $C_1$, ... , $C_k$ translatés de convexes internes. Encore une fois, en pratique, l'inégalité de contrôle sera vérifiée pour tous les conditionnements $A$ mesurables.

\subsection{Espaces vectoriels localement convexes mesurables}

Etant donné que l'on n'a pas de problème de mesurabilité des applications $\mathfrak{m}_\Lambda$ (contrairement au cas i.i.d.), on étend la notion d'espace vectoriel localement convexe mesurable précédemment définie. Soient $X$ un espace vectoriel réel, $\mathcal{F}$ une tribu sur $X$ stable par dilatation et translation. On suppose qu'il existe une famille $\mathcal{C}_0$ de parties de $X$ telle que :

\medskip

\textup{\textsf{(EVLCM$_1$)}} pour tout $C \in \mathcal{C}_0$, $C$ est un convexe interne mesurable (\emph{i.e.} $\mathcal{C}_0 \subset \mathcal{C}_0(\mathcal{F})$) et symétrique ;

\medskip

\textup{\textsf{(EVLCM$_2$)}} $\mathcal{C}_0$ est stable par intersection finie et par dilatation de rapport non nul, et l'intersection des éléments de $\mathcal{C}_0$ est réduite à $\{ 0 \}$ ;

\medskip

Alors, $\mathcal{C}_0$ est un système fondamental de voisinages de $0$ pour une unique topologie localement convexe (séparée) $\tau$. Si, de plus

\medskip

\textup{\textsf{(EVLCM$_3$)}} toute forme linéaire continue sur $X$ est mesurable,

\medskip

on dit que $(X, \mathcal{C}_0, \mathcal{F}, \tau)$ est un \emph{espace vectoriel localement convexe mesurable} (\emph{e.v.l.c.m.}).


\medskip

\textbf{Remarque :} On a supprimé deux hypothèses : on n'impose plus la séparabilité de la topologie localement convexe engendrée par $C \in \mathcal{C}_0$ ; et on ne demande pas que $\mathcal{F}$ soit engendrée par les translatés des éléments de $\mathcal{C}_0$.

\medskip

On définit la notion de limite projective d'e.v.l.c.m. Soient $X$ un espace vectoriel réel, $\mathcal{F}$ une tribu sur $X$ stable par dilatation et translation, et
$$
\overleftarrow{\mathcal{X}} = \big( X_i, f_i, f_{ij} \big)_{i \leqslant j}
$$
une famille telle que

\medskip

\textsf{(PROJ$_1$)} les indices $i$ et $j$ décrivent un ensemble $(J, \leqslant)$ préordonné filtrant à droite ;

\medskip

\textsf{(PROJ$_2$)} pour tout $i \in J$, $(X_i, \mathcal{C}_{0,i}, \mathcal{F}_i, \tau_i)$ est un e.v.l.c.m., $f_i$ une application linéaire mesurable de $X$ dans $X_i$ et, pour tout $(i, j) \in J^2$ tel que $i \leqslant j$, $f_{ij}$ est une application linéaire, continue et mesurable de $X_j$ dans $X_i$ ;

\medskip

\textsf{(PROJ$_3$)} pour tout $(i, j, k) \in J^3$ tel que $i \leqslant j \leqslant k$, on a $f_{ii} = id_{X_i}$,
$$
f_i = f_{ij} \circ f_j  \qquad \textrm{et} \qquad f_{ik} = f_{ij} \circ f_{jk}
$$

On dit alors que $\smash{\overleftarrow{\mathcal{X}}}$ est un \emph{système projectif d'espaces vectoriels localement convexes mesurables}. De plus, si l'on note
$$
\mathcal{C}_0 = \left\{ f_i^{-1}(C_i) \, ; \, i \in J, \, C_i \in \mathcal{C}_{0}(\mathcal{F}_i) \right\}
$$
et $\tau$ la topologie localement convexe engendrée par $\mathcal{C}_0$, $(X, \mathcal{C}_0, \mathcal{F}, \tau)$ est un e.v.l.c.m., appelé \emph{limite projective}\footnote{Au niveau des ensembles, cette définition de limite projective d'ensembles est un peu plus générale que celle de \cite[III.51.]{BouE}} du système projectif d'e.v.l.c.m. $\smash{\overleftarrow{\mathcal{X}}}$.

\section{Théorie de Cramér pour des champs asymptotiquement découplés}

Dans toute cette partie, $(X, \mathcal{C}_0, \mathcal{F}, \tau)$ désigne un e.v.l.c.m., $d$ un entier strictement positif et $(\mu_\Lambda)_{\Lambda \in \mathbb{B}}$ un s.c.i.t. sur $\smash{X^{\mathbb{Z}^d}}$. Pour tout $x \in X$, on note
$$
\mathcal{C}_x = \{ x + C \, ; \, C \in \mathcal{C}_0 \}
$$
et, pour tout $n \geqslant 1$,
$$
\mu_n = \mu_{\Lambda(n)} \circ \big( \mathfrak{m}_{\Lambda(n)} \big)^{-1}
$$
On dira que $(\mu_n)_{n \geqslant 1}$ est la \emph{suite de Cramér} associée au s.c.i.t. $(\mu_\Lambda)_{\Lambda \in \mathbb{B}}$. Si $D$ est une partie convexe de $X$ et si $(\mu_\Lambda)_{\Lambda \in \mathbb{B}}$ est porté par $D^{\mathbb{Z}^d}$, alors, pour tout $n \geqslant 1$, $\mu_n$ est une probabilité portée par $D$, au sens où, pour tout $(A, B) \in \mathcal{F}^2$,
$$
A \cap D = B \cap D \Rightarrow \mu_n(A) = \mu_n(B)
$$

\subsection{Entropie et pression}

On appelle \emph{entropie} de $(\mu_n)_{n \geqslant 1}$ la fonction $s$ définie par : pour tout $x \in X$,
$$
s(x) := \underline{s}(x) := \inf_{C \in \mathcal{C}_x} \liminf_{n \to \infty} \frac{1}{|\Lambda(n)|} \log \mu_n(C)
$$
On définit également l'\emph{entropie supérieure} de $(\mu_n)_{n \geqslant 1}$ par : pour tout $x \in X$,
$$
\overline{s}(x) := \inf_{C \in \mathcal{C}_x} \limsup_{n \to \infty} \frac{1}{|\Lambda(n)|} \log \mu_n(C)
$$
Par construction, l'entropie $s$ est la plus grande fonction vérifiant la borne inférieure :

\medskip

\textsf{(BI)} pour tout $A \in \mathcal{F}$,
$$
\liminf_{n \to \infty} \frac{1}{|\Lambda(n)|} \log \mu_n(A) \geqslant \sup_{\aroi} s
$$

On dit que $(\mu_n)_{n \geqslant 1}$ vérifie un \emph{principe de grandes déviations} (PGD) si la borne supérieure suivante est vérifiée :

\medskip

\textsf{(BS)} pour tout $A \in \mathcal{F}$,
$$
\limsup_{n \to \infty} \frac{1}{|\Lambda(n)|} \log \mu_n(A) \leqslant \sup_{\smash{\overline{A}}} s
$$

De manière générale, \textsf{(BS)} n'est pas vérifiée pour tous les mesurables. Si $\mathcal{P}$ désigne un ensemble de parties mesurables de $X$, on définit la version restreinte de la borne supérieure suivante :

\medskip

\textsf{(BS$_\mathcal{P}$)} pour tout $A \in \mathcal{P}$,
$$
\limsup_{n \to \infty} \frac{1}{|\Lambda(n)|} \log \mu_n(A) \leqslant \sup_{\smash{\overline{A}}} s
$$

En particulier, si $D$ est une partie de $X$, on notera :

\medskip

\textsf{(BS$_{\flat , D}$)} pour tout $K \in \mathcal{F}$  tel que $K \cap D$ soit relativement compact,
$$
\limsup_{n \to \infty} \frac{1}{|\Lambda(n)|} \log \mu_n(K) \leqslant \sup_{\smash{\overline{K}}} s
$$

Si $(\mu_n)_{n \geqslant 1}$ vérifie \textsf{(BS$_{\flat , X}$)}, on dit que $(\mu_n)_{n \geqslant 1}$ vérifie un \emph{principe de grandes déviations faible} (PGD faible). On appelle \emph{pression} de $(\mu_n)_{n \geqslant 1}$ la fonction $p$ définie par : pour tout $\lambda \in X^*$ (mesurable, par hypothèse),
$$
p(\lambda) := \limsup_{n \to \infty} \frac{1}{|\Lambda(n)|} \log \int e^{|\Lambda(n)| \langle \lambda | x \rangle} d\mu_n(x)
$$

\subsection{\'Enoncé des résultats principaux}

Si $(\mu_\Lambda)_{\Lambda \in \mathbb{B}}$ est un s.c.i.t. et si $D$ est une partie de $X$, on dit que $(\mu_\Lambda)_{\Lambda \in \mathbb{B}}$ est \emph{convexe-tendu localement sur $D$} (\emph{c.t.l. sur $D$}) si, pour tout $\gamma > 0$, il existe $K(\gamma) \in \mathcal{F}$ tel que $K(\gamma) \cap D$ soit convexe relativement compact et tel que, pour toute boîte $\Lambda$ contenant $0$,
$$
\mathbb{P} \big( \sigma(0) \in K(\gamma) \big| \mathcal{F}_{\Lambda \setminus \{ 0 \} } \big) \geqslant 1 - \gamma
$$

\begin{Th}
Soient $(X, \mathcal{C}_0, \mathcal{F}, \tau)$ un e.v.l.c.m., $d$ un entier strictement positif, $D$ une partie convexe de $X$, $(\mu_\Lambda)_{\Lambda \in \mathbb{B}}$ un s.c.i.t. sur $\smash{X^{\mathbb{Z}^d}}$ porté par $\smash{D^{\mathbb{Z}^d}}$, et $(\mu_n)_{n \geqslant 1}$ la suite de Cramér associée. Si $(\mu_\Lambda)_{\Lambda \in \mathbb{B}}$ est a.d.i. et c.l., alors\\
$\bullet$ $(\mu_n)_{n \geqslant 1}$ vérifie \textup{\textsf{(BS$_{\flat , D}$)}} ; en particulier, $(\mu_n)_{n \geqslant 1}$ et $(\mu_n|_D)_{n \geqslant 1}$ vérifient un PGD faible ;\\
$\bullet$ la pression de $(\mu_n)_{n \geqslant 1}$ est définie par
$$
\forall \lambda \in X^* \qquad p(\lambda) = \lim_{n \to \infty} \frac{1}{|\Lambda(n)|} \log \int e^{|\Lambda(n)| \langle \lambda | x \rangle} d\mu_n(x)
$$
et vérifie
$$
\forall \lambda \in X^* \quad \forall x \in X \qquad p(\lambda) - s(x) \geqslant \langle \lambda | x \rangle
$$
$\bullet$ si, de plus, $(\mu_\Lambda)_{\Lambda \in \mathbb{B}}$ est c.t.l. sur $D$, alors l'entropie est l'opposée de la convexe-conjuguée de la pression, autrement dit
$$
\forall x \in X \qquad s(x) = \inf_{\lambda \in X^*} \big( p(\lambda) - \langle \lambda | x \rangle \big)
$$
\end{Th}

Si $f$ est une application de $X$ dans un ensemble $E$, on définit
$$
f^{\mathbb{Z}^d} : \big( \sigma(z) \big)_{z \in \mathbb{Z}^d}  \in X^{\mathbb{Z}^d} \mapsto \big( f(\sigma(z)) \big)_{z \in \mathbb{Z}^d} \in E^{\mathbb{Z}^d}
$$
Le théorème précédent s'étend aux limites projectives en utilisant le théorème suivant (dont on trouve une démonstration dans \cite{Pet11c}) :

\begin{Th}[Dawson-Gärtner linéaire]
Soient $\smash{\overleftarrow{\mathcal{X}} = (X_i, f_i, f_{ij})_{(i, j) \in J^2}}$ un système projectif d'espaces vectoriels localement convexes mesurables et $(X, \mathcal{C}_0, \mathcal{F}, \tau)$ sa limite projective. Soient $(\mu_\Lambda)_{\Lambda \in \mathbb{B}}$ un s.c.i.t. sur $X^{\mathbb{Z}^d}$ et $(\mu_n)_{n \geqslant 1}$ la suite de Cramér associée à $(\mu_\Lambda)_{\Lambda \in \mathbb{B}}$. Notant $s$ (resp. $s_i$, $\overline{s}$, $\overline{s}_i$, $p$, $p_i$) l'entropie de $(\mu_n)_{n \geqslant 1}$ (resp. l'entropie de $(\mu_n \circ f_i^{-1})_{n \geqslant 1}$, l'entropie supérieure de $(\mu_n)_{n \geqslant 1}$, l'entropie supérieure de $(\mu_n \circ f_i^{-1})_{n \geqslant 1}$, la pression de $(\mu_n)_{n \geqslant 1}$, la pression de $(\mu_n \circ f_i^{-1})_{n \geqslant 1}$), on a :
$$
s = \inf_{i \in J} (s_i \circ f_i), \quad \overline{s} = \inf_{i \in J} (\overline{s}_i \circ f_i) \quad \textrm{et} \quad - p^* = \inf_{i \in J} (- p_i^* \circ f_i)
$$
En particulier, si, pour tout $i \in J$, $s_i = \overline{s}_i$, alors $(\mu_n)_{n \geqslant 1}$ vérifie \textup{\textsf{(BS$_{\flat , D}$)}}. Si, de plus, pour tout $i \in J$, $s_i = -p_i^*$, alors $s = -p^*$.
\end{Th}

\begin{Th}[Limites projectives]
Soient $\smash{\overleftarrow{\mathcal{X}} = (X_i, f_i, f_{ij})_{(i, j) \in J^2}}$ un système projectif d'e.v.l.c.m. et $(X, \mathcal{C}_0, \mathcal{F}, \tau)$ sa limite projective. Soient $d$ un entier strictement positif, $D$ une partie convexe de $X$, $(\mu_\Lambda)_{\Lambda \in \mathbb{B}}$ un s.c.i.t. sur $\smash{X^{\mathbb{Z}^d}}$ porté par $\smash{D^{\mathbb{Z}^d}}$ et $(\mu_n)_{n \geqslant 1}$ la suite de Cramér associée. Si, pour tout $i \in J$, $( \mu_\Lambda \circ (f_i^{\mathbb{Z}^d})^{-1} )_{\Lambda \in \mathbb{B}}$ est a.d.i. et c.l., alors\\
$\bullet$ $(\mu_n)_{n \geqslant 1}$ vérifie \textup{\textsf{(BS$_{\flat , D}$)}} ; en particulier, $(\mu_n)_{n \geqslant 1}$ et $(\mu_n|_D)_{n \geqslant 1}$ vérifient un PGD faible ;\\
$\bullet$ la pression de $(\mu_n)_{n \geqslant 1}$ est définie par
$$
\forall \lambda \in X^* \qquad p(\lambda) = \lim_{n \to \infty} \frac{1}{|\Lambda(n)|} \log \int e^{|\Lambda(n)| \langle \lambda | x \rangle} d\mu_n(x)
$$
et vérifie
$$
\forall \lambda \in X^* \quad \forall x \in X \qquad p(\lambda) - s(x) \geqslant \langle \lambda | x \rangle
$$
$\bullet$ Si, de plus, pour tout $i \in J$, $( \mu_\Lambda \circ (f_i^{\mathbb{Z}^d})^{-1} )_{\Lambda \in \mathbb{B}}$ est c.t.l. sur $f_i(D)$, alors l'entropie est l'opposée de la convexe-conjuguée de la pression, autrement dit
$$
\forall x \in X \qquad s(x) = \inf_{\lambda \in X^*} \big( p(\lambda) - \langle \lambda | x \rangle \big)
$$
\end{Th}

\subsection{Exemples d'applications}

Le cadre défini ici contient la théorie de Cramér dans le cas indépendant de \cite{Pet11a} et le cadre de \cite{Pfi01}. En particulier, on a donc le PGD faible et l'égalité $s = -p^*$ dans le cadre de la théorie de Cramér dans les espaces de Banach séparables, ainsi qu'en topologie faible, et dans le cadre du théorème de Sanov. D'autre part, pour les champs non indépendants, la théorie développée dans ce chapitre s'applique immédiatement aux exemples mentionnés dans \cite{Pfi01} (champs de Gibbs, chaînes de Markov, \emph{etc.}) pour les principes de grandes déviations de niveau 2 ou 3 (dans la terminologie de \cite{Ell85}), car la relative compacité assure le contrôle local et la convexe-tension locale. Pour les résultats de niveau 1, le contrôle local est vérifié notamment si $\eta$ prend ses valeurs dans un compact ou si $\eta$ est a.d.i. avec $g(1) = 0$.




\section{Retour sur les hypothèses et compléments}

Dans toute cette partie, $(X, \mathcal{C}_0, \mathcal{F}, \tau)$ désigne un e.v.l.c.m., $d$ un entier strictement positif et $(\mu_\Lambda)_{\Lambda \in \mathbb{B}}$ un s.c.i.t. sur $\smash{X^{\mathbb{Z}^d}}$. Pour tout $x \in X$, on note
$$
\mathcal{C}_x = \{ x + C \, ; \, C \in \mathcal{C}_0 \}
$$
et, pour tout $n \geqslant 1$,
$$
\mu_n = \mu_{\Lambda(n)} \circ \big( \mathfrak{m}_{\Lambda(n)} \big)^{-1}
$$
la suite de Cramér associée.

\subsection{Partition adaptée de $\mathbb{Z}^d$}

Nous décrivons, dans cette section, une partition des sites de $\mathbb{Z}^d$ adaptée au découplage asymptotique (nous reprenons la construction de Pfister). Soient $m$ et $n$ deux entiers tels que $n \geqslant m \geqslant 1$. La division euclidienne de $n$ par $m+g(m)$ s'écrit
$$
n = q \big( m+g(m) \big) + r
$$
On pave (pas complètement) $\Lambda(n)$ avec $q^d$ boîtes isométriques à $\Lambda(m+g(m))$, que l'on note $\Lambda'_k$, pour $k \in [1, q^d]$. Puis, pour chaque $k$, on note $\Lambda_k$ la boîte isométrique à $\Lambda(m)$ dont le plus petit élément pour l'ordre lexicographique usuel de $\mathbb{Z}^d$ --- le ``coin en bas à gauche''\,--- coïncide avec le coin en bas à gauche de $\Lambda'_k$. Soit enfin
$$
S_0 = \Lambda(n) \setminus \bigcup_{k=1}^{q^d} \Lambda_k
$$
l'ensemble des sites marginaux. Les boîtes $\Lambda_k$ dépendent de $n$ et $m$ : ces indices seront sous-entendus.

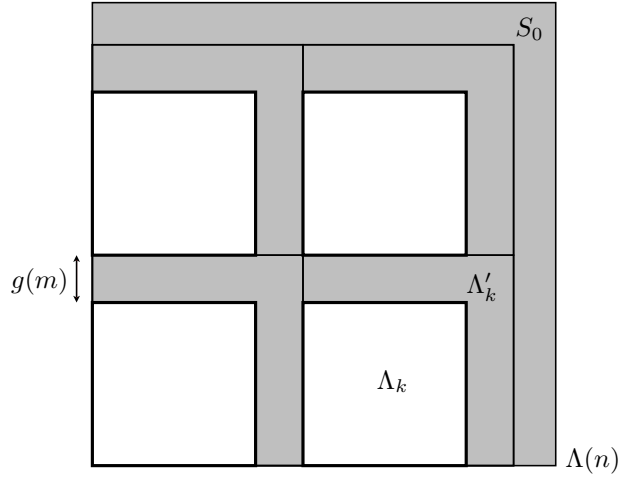
\begin{center}
\def\JPicScale{.7}
\input{lambda2.pst}
\end{center}

On notera que la distance entre deux boîtes $\Lambda_k$ distinctes est supérieure à $g(m)$, ce qui permet d'utiliser le découplage asymptotique et d'énoncer :
\begin{Pro}\label{adifonc}
Pour tout $k \in \{ 1, \ldots , q^d \}$, soit $f_k : X \rightarrow [-\infty , +\infty]$ mesurable. On a :
$$
\mathbb{E} \left( \prod_{k=1}^{q^d} e^{f_k(\mathfrak{m}_{\Lambda_k}(\sigma))} \right) \geqslant e^{-q^{d} c(m)} \prod_{k=1}^{q^d} \mathbb{E} \Big( e^{f_k(\mathfrak{m}_{\Lambda_k}(\sigma))} \Big)
$$
\end{Pro}

\textbf{Démonstration :} On applique $q^d - 1$ fois l'hypothèse de découplage asymptotique inférieur :
\begin{align*}
\phantom\qquad\qquad\quad\mathbb{E} \left( \prod_{k=1}^{q^d} e^{f_k(\mathfrak{m}_{\Lambda_k}(\sigma))} \right)
 &= \int_{t_1, \ldots , t_{q^d} \geqslant 0} \mathbb{P} \Big( \forall k \in \{ 1, \ldots , q^d \} \quad e^{f_k(\mathfrak{m}_{\Lambda_k}(\sigma))} > t_k \Big)\\
 &\geqslant \int_{t_1, \ldots , t_{k^d} \geqslant 0} \big(e^{-c(m)}\big)^{q^d-1} \prod_{k=1}^{q^d} \mathbb{P}\Big( e^{f_k(\mathfrak{m}_{\Lambda_k}(\sigma))} > t_k \Big)\\
 &\geqslant e^{-q^{d} c(m)} \prod_{k=1}^{q^d} \mathbb{E} \Big( e^{f_k(\mathfrak{m}_{\Lambda_k}(\sigma))} \Big)\qquad\qquad\qquad\qquad\qquad\,\qed
\end{align*}

\begin{Pro}\label{clfonc}
Soit $f : X \rightarrow [-\infty , +\infty]$ une application mesurable. On suppose qu'il existe $\beta \in \mathbb{R}$ et $C \in \mathcal{C}_0(\mathcal{F})$ tels que
$$
\{ x \in X \, ; \, f(x) \geqslant \beta \} \supset t(C) \cdot C
$$
Alors, en notant $S = \Lambda_1 \cup \cdots \cup \Lambda_{q^d}$,
$$
\mathbb{E} \left( \prod_{z \in S_0} e^{f(\sigma(z))} \bigg| \mathcal{F}_S \right) \geqslant \big( e^\beta \alpha(C) \big)^{|S_0|}
$$
Si, de plus, $(\mu_\Lambda)_{\Lambda \in \mathbb{B}}$ est porté par $D^{\mathbb{Z}^d}$ et c.t.l. sur $D$, et si $\gamma \leqslant \alpha(C)/2$, alors
$$
\mathbb{E} \left( \prod_{z \in S_0} e^{f(\sigma(z))} 1_{ \sigma(z) \in K(\gamma) } \bigg| \mathcal{F}_S \right) \geqslant \big( e^\beta \alpha(C)/2 \big)^{|S_0|}
$$
\end{Pro}

\textbf{Démonstration :} Soient $z_0 \in S_0$ et $S_1 = S \cup (S_0 \setminus \{ z_0 \})$. On a :
$$
\mathbb{E} \left( \prod_{z \in S_0} e^{f(\sigma(z))} \bigg| \mathcal{F}_S \right) = \mathbb{E} \left( \mathbb{E} \big( e^{f(\sigma(z_0))} \big| \mathcal{F}_{S_1} \big) \prod_{z \in S_0 \setminus \{ z_0 \}} e^{f(\sigma(z))} \bigg| \mathcal{F}_S \right)
$$
Or, l'inégalité de Tchebychev donne
$$
\mathbb{E} \big( e^{f(\sigma(z_0))} \big| \mathcal{F}_{S_1} \big) \geqslant e^\beta \mathbb{P} \big( f(\sigma(z_0)) \geqslant \beta \big| \mathcal{F}_{S_1} \big) \geqslant e^\beta \mathbb{P} \big( \sigma(z_0) \in t(C) \cdot C \big| \mathcal{F}_{S_1} \big) \geqslant e^\beta \alpha(C)
$$
Poursuivant par récurrence, on obtient le premier résultat. Quant au second, il suffit d'adapter la preuve précédente ainsi :
$$
\mathbb{E} \big( e^{f(\sigma(z_0))} 1_{ \sigma(z_0) \in K(\gamma) } \big| \mathcal{F}_{S_1} \big) \geqslant e^\beta \mathbb{P} \big( \sigma(z_0) \in \big( t(C) \cdot C \big) \cap K(\gamma) \big| \mathcal{F}_{S_1} \big) \geqslant e^\beta \alpha(C)/2
$$
et de continuer par récurrence.\qed

\medskip

Notons
$$
\rho_{m, n} := \frac{|S_0|}{|\Lambda(n)|}
$$
la proportion de sites marginaux dans $\Lambda(n)$. Un petit calcul donne
$$
\rho_{m, n} = 1 - \left( 1 - \frac{r}{n} \right)^d \left( 1 - \frac{g(m)}{m + g(m)} \right)^d \lesssim d \cdot \left( \frac{m}{n} + \frac{g(m)}{m + g(m)} \right) \xrightarrow[\substack{m \to \infty\\n \geqslant m}]{} 0
$$
Remarquons que, pour que $\rho_{m, n}$ puisse être rendu arbitrairement petit, il faut que le second terme de la dernière parenthèse tende vers $0$ avec $m$, ce qui équivaut à
$$
\frac{g(m)}{m} \to 0
$$
L'hypothèse faite sur $g$ permet donc que la proportion de sites marginaux soit asymptotiquement négligeable. Ainsi l'hypothèse de contrôle local s'avérera suffisante pour contrôler les sites marginaux.

\subsection{Convexe-tension locale}

On rappelle que, si $(\mu_\Lambda)_{\Lambda \in \mathbb{B}}$ est un s.c.i.t. et si $D$ est une partie de $X$, on dit que $(\mu_\Lambda)_{\Lambda \in \mathbb{B}}$ est \emph{convexe-tendu localement sur $D$} (\emph{c.t.l. sur $D$}) si, pour tout $\gamma > 0$, il existe $K(\gamma) \in \mathcal{F}$ tel que $K(\gamma) \cap D$ soit convexe relativement compact et tel que, pour toute boîte $\Lambda$ contenant $0$,
$$
\mathbb{P} \big( \sigma(0) \in K(\gamma) \big| \mathcal{F}_{\Lambda \setminus \{ 0 \} } \big) \geqslant 1 - \gamma
$$
En particulier, passant aux espérances, on voit que la loi $\mu_1$ de $\sigma(0)$ est convexe-tendue sur $D$ (au sens entendu dans le cas i.i.d.). Plus généralement :
\begin{Pro}\label{ct}
Si $(\mu_\Lambda)_{\Lambda \in \mathbb{B}}$ est c.t.l. sur $D$, alors, pour tout $n \geqslant 1$, $\mu_n$ est convexe-tendue sur $D$. Plus précisément,
$$
\lim_{\gamma \to 0^+} \mu_n \big( K(\gamma) \big) = 1
$$
\end{Pro}

\textbf{Démonstration :} Soient $n \in \mathbb{N}$ et $\gamma > 0$. Soit $K = K(\gamma/n^d)$. Alors
$$
\phantom\qquad\qquad \mu_n(X \setminus K) \leqslant \mu_{\Lambda(n)} \big( \exists z \in \Lambda(n) \quad \sigma(z) \in X \setminus K \big) \leqslant n^d \mu_1(X \setminus K) < \gamma \qquad\quad\qed
$$

\section{Démonstration des résultats principaux}

On garde les notations introduites dans la partie précédente. Sachant que, dans la partie \'Enoncé des résultats principaux, le théorème sur les limites projectives est conséquence du premier et du théorème de Dawson-Gärtner linéaire, on donne uniquement la démonstration du premier.

\subsection{Entropie}

\begin{Th}[lemme sous-additif]\label{lsa}
Si $(\mu_\Lambda)_{\Lambda \in \mathbb{B}}$ est a.d.i. et c.l., alors, pour tout $x \in X$, pour tout $C \in \mathcal{F}$ convexe contenant un convexe $B \in \mathcal{C}_0(\mathcal{F})$ et pour tout $\varepsilon > 0$, il existe $M(\varepsilon) \in \mathbb{N}$ tel que, pour tous $n \geqslant m \geqslant M(\varepsilon)$,
$$
\frac{1}{|\Lambda(n)|} \log \mu_n \big( x + (1+\varepsilon)C \big) \geqslant \frac{1}{|\Lambda(m)|} \log \mu_m(x + C) - \frac{c(m)}{\Lambda(m)} + \rho_{m, n} \log \alpha(B)
$$
\end{Th}


\textbf{Démonstration :} Quitte à travailler sur les champs $(\sigma(z) - x)_{z \in \mathbb{Z}^d}$, on se ramène à $x = 0$. Soient $C \in \mathcal{F}$ convexe contenant un convexe $B \in \mathcal{C}_0(\mathcal{F})$ et $m, n \in \mathbb{N}^*$ avec $m \leqslant n$. La sous-additivité et la $\mathbb{R}_+$-linéarité de $M_C$ montrent que
\begin{align*}
M_C \big( \mathfrak{m}_{\Lambda(n)} (\sigma) \big) &\leqslant \sum_{z \in S_0} \frac{1}{|\Lambda(n)|} M_C \big( \sigma(z) \big) + \sum_{k=1}^{q^d} \frac{|\Lambda(m)|}{|\Lambda(n)|} M_C \big( \mathfrak{m}_{\Lambda_k} (\sigma) \big)\\
 &\leqslant \sum_{z \in S_0} \frac{1}{|S_0|} \rho_{m, n} M_C \big( \sigma(z) \big) + \sum_{k=1}^{q^d} \frac{1}{q^d} M_C \big( \mathfrak{m}_{\Lambda_k} (\sigma) \big)
\end{align*}
D'où, pour $\varepsilon > 0$,
\begin{align*}
\frac{1}{|\Lambda(n)|}& \log \mu_n\big( (1+\varepsilon)C \big) \\
 &\geqslant \frac{1}{|\Lambda(n)|} \log \mathbb{P}\Big(M_C \big( \mathfrak{m}_{\Lambda(n)} (\sigma) \big) < 1+\varepsilon \Big)\\
 &\geqslant \frac{1}{|\Lambda(n)|} \log \mathbb{P}\Big(\forall z \in S_0 \quad \rho_{m, n} M_C \big(\sigma(z)\big) < \varepsilon \, ; \\
 & \qquad\qquad\qquad\qquad \forall k \in [1, q^d] \quad M_C \big(\mathfrak{m}_{\Lambda_k} (\sigma) \big) \leqslant 1 \Big)\\
 &\geqslant \frac{1}{|\Lambda(n)|} \log \mathbb{E} \left( \mathbb{E} \bigg( \prod_{z \in S_0} \exp \big( \delta_{B_1} (\sigma(z)) \big) \bigg| \mathcal{F}_S \bigg) \prod_{k=1}^{q^d} \exp \big( \delta_C (\mathfrak{m}_{\Lambda_k} (\sigma) ) \big) \right)
\end{align*}
où on a noté $S = \Lambda_1 \cup \cdots \Lambda_{q^d}$, $B_1 = (\varepsilon / \rho_{m, n}) B$ et $\delta_{B_1} = -\infty 1_{X \setminus B_1}$. Compte tenu de la majoration de $\rho_{m, n}$ donnée dans la section précédente, il existe $M(\varepsilon) \in \mathbb{N}$ tel que
$$
\forall (m, n) \in \mathbb{N}^2 \qquad n \geqslant m \geqslant M(\varepsilon) \Rightarrow \frac{\varepsilon}{\rho_{m, n}} \geqslant t(B)
$$
Alors, pour $n \geqslant m \geqslant M(\varepsilon)$, la proposition \ref{clfonc}, avec $\beta = 0$, permet de minorer l'espérance conditionnelle et d'obtenir
$$
\frac{1}{|\Lambda(n)|} \log \mu_n \big( (1+\varepsilon)C \big) \geqslant \rho_{m, n} \log \alpha(B) + \frac{1}{|\Lambda(n)|} \log \mathbb{E} \left( \prod_{k=1}^{q^d} \exp \big( \delta_C (\mathfrak{m}_{\Lambda_k} (\sigma) ) \big) \right)
$$
Puis la proposition \ref{adifonc} donne
\begin{align*}
\frac{1}{|\Lambda(n)|} \log \mu_n \big( (1+\varepsilon)C \big) &\geqslant \rho_{m, n} \log \alpha(B) - \frac{q^d c(m)}{|\Lambda(n)|} + \frac{1}{|\Lambda(n)|} \log \prod_{k=1}^{q^d} \mathbb{E} \Big( \exp \big( \delta_C (\mathfrak{m}_{\Lambda_k} (\sigma) ) \big) \Big)\\
 &\geqslant \rho_{m, n} \log \alpha(B) - \frac{c(m)}{|\Lambda(m)|} + \frac{1}{|\Lambda(m)|} \log \mu_m(C)
\end{align*}
où on a utilisé le fait que
$$
\frac{q^d}{|\Lambda(n)|} = \frac{1 - \rho_{m, n}}{|\Lambda(m)|} \leqslant \frac{1}{|\Lambda(m)|}
$$
et l'invariance par translation.\qed

\medskip

\textbf{Remarque :} Intuitivement, on a d'abord fait en sorte que l'écartement entre les petites boîtes soit suffisamment petit devant la taille de ces boîtes.
Puis, on a pris suffisamment de petites boîtes pour que la proportion occupée par elles dans une grande boîte soit assez grande. Au final, on obtient une sous-additivité asymptotique entre les petites boîtes et la grande.

\begin{Cor}\label{lsacor}
Si $(\mu_\Lambda)_{\Lambda \in \mathbb{B}}$ est a.d.i. et c.l., on a l'égalité $\underline{s} = \overline{s}$.
\end{Cor}

\textbf{Démonstration :} Prenant, dans l'inégalité du théorème précédent, et pour $C \in \mathcal{C}_0$, la limite inférieure en $n$, puis la limite supérieure en $m$, on obtient
$$
\liminf_{n \to \infty} \frac{1}{|\Lambda(n)|} \log \mu_n \big(x + (1+\varepsilon)C \big) \geqslant \limsup_{m \to \infty} \frac{1}{|\Lambda(m)|} \log \mu_m(x + C)
$$
Reste à prendre l'infimum en $C \in \mathcal{C}_0$ pour conclure.\qed

\begin{Pro}\label{sconc}
L'entropie $s$ est semi-continue supérieurement et concave.
\end{Pro}
\textbf{Démonstration :} Elle est analogue à celle dans le cas i.i.d. Montrons la semi-continuité supérieure. Soient $t \in \mathbb{R}$ et $x \in X$ tels que $s(x) < t$. Il existe $C \in \mathcal{C}_x$ tel que
$$
\liminf_{n \to \infty} \frac{1}{|\Lambda(n)|} \log \mu_n(C) < t
$$
Alors, pour tout $y \in C$,
$$
s(y) \leqslant \liminf_{n \to \infty} \frac{1}{|\Lambda(n)|} \log \mu_n(C) < t
$$
donc $x \in C \subset \{ s < t \}$ et $\{ s < t \}$ est ouvert. Pour la concavité, il suffit alors de prouver que
$$
s\left( \frac{x+y}{2} \right) \geqslant \frac{s(x)+s(y)}{2}
$$
La semi-continuité supérieure permet alors de passer des inégalités
$$
s \big( (1 - u)x + uy \big) \geqslant (1 - u)s(x)+ us(y)
$$
pour $u$ dyadique à celles pour $u \in [0, 1]$. La preuve est un raffinement de la démonstration du lemme sous-additif (\ref{lsa}). Soit $C \in \mathcal{C}_0$. Par continuité des applications d'espace vectoriel, pour tout $\varepsilon \in ]0, 1[$, il existe $C_x \in \mathcal{C}_x$ et $C_y \in \mathcal{C}_y$ tels que
$$
\frac{1}{2}C_x + \frac{1}{2}C_y \subset \frac{1}{2}(x+y) + (1-\varepsilon)C
$$
Puis
\begin{align*}
\frac{1}{|\Lambda(n)|} \log &\, \mu_n \left( \frac{1}{2}(x+y) + C \right) \\
 &\geqslant \frac{1}{|\Lambda(n)|} \log \mathbb{P} \Bigg(\forall z \in S_0 \quad \rho_{m, n} M_C(\sigma(z)) < \varepsilon \, ;\\
&\qquad\qquad\qquad\qquad\bigcap_{q \textrm{ pair}} \{ \mathfrak{m}_{\Lambda_q} \sigma \in C_x \} \cap \bigcap_{q \textrm{ impair}} \{ \mathfrak{m}_{\Lambda_q} \sigma \in C_y \} \Bigg)
\end{align*}
En poursuivant comme précédemment, on obtient
\begin{align*}
\liminf_{n \to \infty} \frac{1}{|\Lambda(n)|} & \log \mu_n \left( \frac{1}{2}(x+y) + C \right) \\
 & \geqslant \frac{1}{2} \left( \liminf_{n \to \infty} \frac{1}{|\Lambda(n)|} \log \mu_n(C_x) + \liminf_{n \to \infty} \frac{1}{|\Lambda(n)|} \log \mu_n(C_y) \right) \\
 & \geqslant \frac{1}{2} \big( s(x) + s(y) \big)
\end{align*}
puis l'inégalité désirée en passant à l'infimum en $C \in \mathcal{C}_0$.\qed

\medskip

Pour achever la démonstration du PGD faible, on montre un résultat un peu plus fort. On étend l'addition à $[-\infty , +\infty]$ via
$$
-\infty \,\dot{+}\, (+\infty) = +\infty
$$
Si $f$ est une fonction de $X$ dans $[-\infty , +\infty]$, on dit que $f$ est \emph{semi-continue supérieurement} si, pour tout $t \in [-\infty , +\infty]$, l'ensemble
$$
\big\{ x \in X \, ; \, f(x) < t \big\}
$$
est un ouvert de $X$. Pour toute fonction $f : X \rightarrow [-\infty , +\infty]$, on note $f^\bullet$ sa régularisée semi-continue supérieurement, autrement dit la plus petite fonction semi-continue supérieurement qui soit supérieure à $f$.
\begin{Th}[Varadhan compact]\label{varadcompadi}
Soit $D$ une partie convexe de $X$. On suppose $(\mu_\Lambda)_{\Lambda \in \mathbb{B}}$ a.d.i., c.l. et porté par $D^{\mathbb{Z}^d}$. Pour toute fonction mesurable $f : X \rightarrow [-\infty , +\infty]$ et pour tout $K \in \mathcal{F}$ tel que $K \cap D$ soit relativement compact,
$$
\limsup_{n \to \infty} \frac{1}{|\Lambda(n)|} \log \mathbb{E} \Big( e^{|\Lambda(n)| f( \mathfrak{m}_{\Lambda(n)} (\sigma) )} 1_{\mathfrak{m}_{\Lambda(n)} (\sigma) \in K} \Big) \leqslant \sup_{\overline{K \cap D}} [f^\bullet \,\dot{+}\, s]
$$
\end{Th}

\textbf{Démonstration :} Soient $K \in \mathcal{F}$ relativement compact, $\delta > 0$ et $M < 0$. Par définition de $f^\bullet$, pour tout $x \in X$, il existe $C(x) \in \mathcal{C}_x$ tel que, pour tout $y \in C(x)$,
$$
f(y) \leqslant \max \big( f^\bullet(x) + \delta , M \big)
$$
Par définition de $s(x) = \overline{s}(x)$, quitte à réduire un peu $C(x)$, on peut supposer que
$$
\limsup_{n \to \infty} \frac{1}{n} \log \mu_n \big( C(x) \big) \leqslant \max \big( s(x) + \delta , M \big)
$$
Du recouvrement de $\smash{\overline{K \cap D}}$ par les $C(x)$ avec $x \in \smash{\overline{K \cap D}}$, on peut extraire un sous-recouvrement fini, noté $\{ C(x_i) ; i \in \{ 1, \ldots , r \} \}$. Pour tout $n \geqslant 1$, ayant choisi $\sigma$ à valeurs dans $D^{\mathbb{Z}^d}$, on a :
\begin{align*}
\frac{1}{|\Lambda(n)|} & \log \mathbb{E} \Big( e^{|\Lambda(n)| f( \mathfrak{m}_{\Lambda(n)} (\sigma) )} 1_{\mathfrak{m}_{\Lambda(n)} (\sigma) \in K} \Big) \\
 & \leqslant \frac{1}{|\Lambda(n)|} \log \sum_{i=1}^r \mathbb{E} \Big( e^{|\Lambda(n)| f( \mathfrak{m}_{\Lambda(n)} (\sigma) )} 1_{\mathfrak{m}_{\Lambda(n)} (\sigma) \in C(x_i)} \Big)\\
 & \leqslant \frac{1}{|\Lambda(n)|} \log \sum_{i=1}^r e^{|\Lambda(n)| \max \big( f^\bullet(x_i) + \delta, M \big)} \mu_n \big( C(x_i) \big)
\end{align*}
Prenant la limite supérieure en $n$ et utilisant le lemme \ref{ls}, on obtient :
\begin{align*}
\limsup_{n \to \infty} \frac{1}{|\Lambda(n)|} & \log \mathbb{E} \Big( e^{|\Lambda(n)| f( \mathfrak{m}_{\Lambda(n)} (\sigma) )} 1_{\mathfrak{m}_{\Lambda(n)} (\sigma) \in K} \Big) \\
 & \leqslant \max_{1\leqslant i\leqslant r} \Big( \max \big( f^\bullet(x_i) + \delta , M \big) + \max \big( (s(x_i) + \delta) , M \big) \Big)\\
 & \leqslant \sup_{x \in \smash{\overline{K \cap D}}} \Big( \max \big( f^\bullet(x) + \delta , M \big) + \max \big( (s(x) + \delta) , M \big) \Big)
\end{align*}
ce qui donne le résultat attendu, en faisant tendre $\delta$ vers $0$ et $M$ vers $-\infty$.\qed

\begin{Cor}[Borne supérieure faible]
Pour tout $K \in \mathcal{F}$ tel que $K \cap D$ soit relativement compact,
$$
\limsup_{n \to \infty} \frac{1}{|\Lambda(n)|} \log \mu_n (K) \leqslant \sup_{\overline{K \cap D}} s
$$
En particulier, $(\mu_n)_{n \geqslant 1}$ vérifie un PGD faible.
\end{Cor}

\textbf{Démonstration :} Il suffit d'appliquer le lemme de Varadhan compact \ref{varadcompadi} à $f = 0$.\qed

\medskip

Terminons cette section avec un résultat proche de la borne supérieure pour les convexes dans le cas où $(\mu_\Lambda)_{\Lambda \in \mathbb{B}}$ est convexe-tendu localement.

\begin{Th}
Soit $D$ une partie convexe de $X$. On suppose que $(\mu_\Lambda)_{\Lambda \in \mathbb{B}}$ est a.d.i., c.l., porté par $D$ et c.t.l. sur $D$. Alors, pour tout $C \in \mathcal{F}$ convexe contenant un convexe $B \in \mathcal{C}_0(\mathcal{F})$, on a :
$$
\limsup_{n \to \infty} \frac{1}{|\Lambda(n)|} \log \mu_n (C) \leqslant \inf_{\varepsilon > 0} \sup_{(1+\varepsilon)\overline{C}} s
$$
\end{Th}

\textbf{Remarque :} On sait conclure qu'on a bien \textsf{(BS$_c$)} dans le cas où $s$ a ses ensembles de niveau compacts.

\medskip

\textbf{Démonstration :} Soit $\gamma \leqslant \alpha(B)/2$. On adapte la démonstration du lemme sous-additif \ref{lsa} en appliquant la seconde partie de la proposition \ref{clfonc} et on obtient, pour tous $n \geqslant m \geqslant M(\varepsilon)$,
$$
\frac{1}{|\Lambda(m)|} \log \mu_m (K(\gamma) \cap C) - \frac{c(m)}{|\Lambda(m)|} + \rho_{m, n} \log \frac{\alpha(B)}{2} \leqslant \frac{1}{|\Lambda(n)|} \log \mu_n \big( K(\gamma) \cap (1+\varepsilon)C \big)
$$
Puis, la borne supérieure faible donne
$$
\limsup_{n \to \infty} \frac{1}{|\Lambda(n)|} \log \mu_n\big( K(\gamma) \cap (1+\varepsilon)C \big) \leqslant \sup_{\overline{ K(\gamma) \cap (1+\varepsilon)C }} s \leqslant \sup_{ (1+\varepsilon) \overline{C}} s
$$
car $\overline{(1+\varepsilon) C} = (1+\varepsilon) \overline{C}$. Combinant les deux inégalités et prenant le supremum en $\gamma \leqslant \alpha(B)/2$, il vient, sachant que $\mu_m(K(\gamma)) \to 1$ lorsque $\gamma \to 0^+$ (cf. \ref{ct}),
$$
\frac{1}{|\Lambda(m)|} \log \mu_m (C) - \frac{c(m)}{|\Lambda(m)|} + \limsup_{n \to \infty} \rho_{m, n} \log \frac{\alpha(B)}{2} \leqslant \sup_{ (1+\varepsilon) \overline{C}} s
$$
Prenant alors la limite supérieure en $m$, puis l'infimum en $\varepsilon > 0$, on obtient la borne supérieure voulue.\qed

\subsection{Pression}

On rappelle que la pression de la suite $(\mu_n)_{n \geqslant 1}$ est définie par : pour tout $\lambda \in X^*$,
$$
p(\lambda) := \limsup_{n \to \infty} \frac{1}{|\Lambda(n)|} \log \int e^{|\Lambda(n)| \langle \lambda | x \rangle} d\mu_n(x)
$$

\begin{Th}\label{press}
Si $(\mu_\Lambda)_{\Lambda \in \mathbb{B}}$ est a.d.i. et c.l., on a
$$
p(\lambda) = \lim_{n \to \infty} \frac{1}{|\Lambda(n)|} \log \int e^{|\Lambda(n)| \langle \lambda | x \rangle} d\mu_n(x)
$$
\end{Th}

\textbf{Démonstration :} De même que pour l'entropie, nous allons montrer une inégalité sous-additive. Soient $m, n \in \mathbb{N}$ avec $m \leqslant n$. On écrit
\begin{align*}
\frac{1}{|\Lambda(n)|} &\log \mathbb{E} \Big( \exp\big( |\Lambda(n)| \langle \lambda | \mathfrak{m}_{\Lambda(n)} (\sigma) \rangle \big) \Big) \\
 &= \frac{1}{|\Lambda(n)|} \log \mathbb{E} \left( \mathbb{E} \bigg( \prod_{z \in S_0} \exp \big( \langle \lambda | \sigma(z) \rangle \big) \bigg| \mathcal{F}_S \bigg) \prod_{k=1}^{q^d} \exp \big( |\Lambda(m)| \langle \lambda | \mathfrak{m}_{\Lambda_k} (\sigma) \rangle \big) \right)
\end{align*}
avec $S_0 = \Lambda_1 \cup \cdots \cup \Lambda_{q^d}$. Etant donné que $\lambda \in X^*$, il existe $C \in \mathcal{C}_0$ tel que $C \subset \{ x \in X \, ; \, \langle \lambda | x \rangle \geqslant - 1 \}$. On peut donc utiliser, comme pour l'entropie, les propositions \ref{clfonc} (avec $\beta = - t(C)$) puis \ref{adifonc} pour obtenir
\begin{align*}
\frac{1}{|\Lambda(n)|} &\log \mathbb{E} \Big( e^{|\Lambda(n)| \langle \lambda | \mathfrak{m}_{\Lambda(n)} (\sigma) \rangle } \Big) \\
 &\geqslant \rho_{m, n} \big( \log \alpha(C) - t(C) \big) - \frac{c(m)}{|\Lambda(m)|} + (1 - \rho_{m, n}) \frac{1}{|\Lambda(m)|} \log \mathbb{E} \Big( e^{|\Lambda(m)| \langle \lambda | \mathfrak{m}_{\Lambda(m)} (\sigma) \rangle } \Big)
\end{align*}
On conclut alors en prenant la limite inférieure en $n$, puis la limite supérieure en $m$.\qed

\subsection{Pression et entropie}

\begin{Th}\label{spineg}
Pour tout $x \in X$ et pour tout $\lambda \in X^*$, on a :
$$
p(\lambda) - s(x) \geqslant \langle \lambda | x \rangle
$$
\end{Th}

\textbf{Démonstration :} Soient $x \in X$ et $\lambda \in X^*$. Définissons, pour $\varepsilon > 0$, le demi-espace ouvert (mesurable)
$$
H = \{ y \in X ; \langle \lambda | y \rangle > \langle \lambda | x \rangle - \varepsilon \}
$$
L'inégalité de Tchebychev donne
\begin{align*}
\liminf_{n \to \infty} \frac{1}{|\Lambda(n)|} & \log \mu_n(H) \\
 & = \liminf_{n \to \infty}\frac{1}{|\Lambda(n)|} \log \mu_n \big( |\Lambda(n)| \lambda - |\Lambda(n)| \langle \lambda | x \rangle + |\Lambda(n)| \varepsilon > 0 \big)\\
 & \leqslant \liminf_{n \to \infty} \frac{1}{|\Lambda(n)|} \log \int \exp \big( |\Lambda(n)| \lambda - |\Lambda(n)| \langle \lambda | x \rangle + |\Lambda(n)| \varepsilon \big) d\mu_n \\
 & = p(\lambda) - \langle \lambda | x \rangle + \varepsilon
\end{align*}
Etant donné que $H$ est un voisinage mesurable de $x$, il vient $s(x) \leqslant p(\lambda) - \langle \lambda | x \rangle + \varepsilon$. Puis, passant à l'infimum en $\varepsilon > 0$, on obtient
$$
\phantom\qquad\qquad\qquad\qquad\qquad\qquad\qquad\, p(\lambda) - s(x) \geqslant \langle \lambda | x \rangle \qquad\qquad\qquad\qquad\qquad\qquad\qquad\,\qed
$$

Dans le cas où $(\mu_\Lambda)_{\Lambda \in \mathbb{B}}$ est convexe-tendu localement, nous allons montrer que $s = -p^*$.



\begin{Le}
Soit $D$ une partie convexe de $X$. On suppose que $(\mu_\Lambda)_{\Lambda \in \mathbb{B}}$ est a.d.i., c.l., porté par $D$ et c.t.l. sur $D$. Alors la pression est la fonction convexe-conjuguée de l'opposée de l'entropie, \emph{i.e.}
$$
\forall \lambda \in X^* \qquad p(\lambda) = \sup_{x \in X} \big( \langle \lambda | x \rangle + s(x) \big)
$$
\end{Le}

\textbf{Démonstration :} L'inégalité $p \geqslant (-s)^*$ ne nécessite pas la convexe-tension : elle repose sur l'inégalité de Tchebychev. Il suffit de passer au supremum en $x \in X$ dans \ref{spineg} :
$$
\forall \lambda \in X^* \qquad p(\lambda) \geqslant \sup_{x \in X} \big( \langle \lambda | x \rangle + s(x) \big)
$$
Montrons l'autre inégalité. Soient $\lambda \in X^*$ et $m, n \in \mathbb{N}$ avec $m \leqslant n$. Soit $C \in \mathcal{C}_0$ tel que $C \subset \{ x \in X \, ; \, \langle \lambda | x \rangle \geqslant - 1 \}$. Soit $\gamma \leqslant \alpha(C)/2$. On adapte la démonstration du théorème \ref{press} en utilisant la seconde partie de la proposition \ref{clfonc} et on obtient
\begin{align*}
\frac{1}{|\Lambda(m)|} &\log \mathbb{E} \Big( e^{|\Lambda(m)| \langle \lambda | \mathfrak{m}_{\Lambda(m)} (\sigma) \rangle} 1_{\mathfrak{m}_{\Lambda(m)} (\sigma) \in K(\gamma)} \Big) - \frac{c(m)}{|\Lambda(m)|} + \rho_{m, n} \left( \log \frac{\alpha(C)}{2} - t(C) \right) \\
 &\leqslant \frac{1}{|\Lambda(n)|} \log \mathbb{E} \Big( e^{|\Lambda(n)| \langle \lambda | \mathfrak{m}_{\Lambda(n)} (\sigma) \rangle} 1_{\mathfrak{m}_{\Lambda(m)} (\sigma) \in K(\gamma)} \Big)
\end{align*}
Puis, le lemme de Varadhan compact \ref{varadcompadi} donne
$$
\limsup_{n \to \infty} \frac{1}{|\Lambda(n)|} \log \mathbb{E} \Big( e^{|\Lambda(n)| \langle \lambda | \mathfrak{m}_{\Lambda(n)} (\sigma) \rangle} 1_{\mathfrak{m}_{\Lambda(m)} (\sigma) \in K(\gamma)} \Big) \leqslant \sup_{x \in X} \big( \langle \lambda | x \rangle + s(x) \big)
$$
Combinant les deux dernières inégalités et passant au supremum en $\gamma \leqslant \alpha(C)/2$, il vient, sachant que $\mu_m(K(\gamma)) \to 1$ lorsque $\gamma \to 0^+$ (cf. \ref{ct}) et par un argument d'uniforme intégrabilité,
$$
\frac{1}{|\Lambda(m)|} \log \mathbb{E} \Big( e^{|\Lambda(m)| \langle \lambda | \mathfrak{m}_{\Lambda(m)} (\sigma) \rangle} \Big) - \frac{c(m)}{|\Lambda(m)|} + \limsup_{n \to \infty} \rho_{m, n} \left( \log \frac{\alpha(C)}{2} - t(C) \right) \leqslant \sup_{X} (\lambda + s)
$$
Reste alors à prendre la limite en $m$.\qed

\begin{Th}
Soit $D$ une partie convexe de $X$. On suppose que $(\mu_\Lambda)_{\Lambda \in \mathbb{B}}$ est a.d.i., c.l., porté par $D$ et c.t.l. sur $D$. Alors l'entropie est l'opposée de la fonction convexe-conjuguée de la pression, \emph{i.e.}
$$
\forall x \in X \qquad s(x) = \inf_{\lambda \in X^*} \big( p(\lambda) - \langle \lambda | x \rangle \big)
$$
\end{Th}

\textbf{Démonstration :} En reprenant les notations du lemme \ref{fl}, le théorème précédent s'écrit $p = (-s)^*$. Comme $-s$ est convexe et s.c.i. (cf. la proposition \ref{sconc}), donc $\sigma(X, X^*)$-s.c.i., le lemme \ref{fl} donne
$$
\phantom\qquad\qquad\qquad\qquad\qquad\qquad\qquad\, p^* = (-s)^{**} = -s \qquad\qquad\qquad\qquad\qquad\qquad\qquad\,\qed
$$

\newpage

\section{Annexe : les trois fondements de la théorie de Cramér}

Dans cette section, nous rappelons trois résultats généraux sur lesquels reposent la démonstration du cas i.i.d., résultats qui s'étendent au cas asymptotiquement découplé (seul le lemme sous-additif doit être adapté).

\subsection{Lemme sous-additif}

Le premier résultat de ce type remonte à un article de M. Fekete \cite{Fek39}.
\begin{Le}\label{lsa1}
Soit $\big(u(n)\big)_{n\geqslant 1}$ une suite à valeurs dans $[0, +\infty]$. On suppose que

\medskip

\textup{\textsf{(SA)}} $u$ est sous-additive, \emph{i.e.}
$$
\forall m, n \geqslant 1 \qquad u(m+n) \leqslant u(m) + u(n)
$$

Alors,
$$
\liminf_{n\to\infty} \frac{u(n)}{n} = \inf_{n \geqslant 1} \frac{u(n)}{n}
$$
\end{Le}

\textbf{Démonstration :} Par sous-additivité, pour tous $d, m \geqslant 1$,
$$
\frac{u(dm)}{dm} \leqslant \frac{u(m)}{m}
$$
En faisant tendre $d$ vers $\infty$, on obtient
$$
\liminf_{n\to\infty} \frac{u(n)}{n} \leqslant \liminf_{d\to\infty} \frac{u(dm)}{dm} \leqslant \frac{u(m)}{m}
$$
d'où le résultat en passant à l'infimum en $m \geqslant 1$.\qed

\begin{Le}\label{lsa2}
Soit $\big(u(n)\big)_{n\geqslant 1}$ une suite à valeurs dans $[0, +\infty]$. On suppose que

\medskip

\textup{\textsf{(SA)}} $u$ est sous-additive ;

\medskip

\textup{\textsf{(C)}} $u$ est contrôlée, \emph{i.e.} il existe $N \geqslant 1$ tel que
$$
\forall n \geqslant N \qquad u(n) < + \infty
$$

Alors, la suite $\big( u(n)/n \big)_{n\geqslant 1}$ converge vers
$$
\inf_{n \geqslant 1} \frac{u(n)}{n}
$$
\end{Le}

\textbf{Démonstration :} Soient $n \geqslant m \geqslant N$. La division euclidienne de $n$ par $m$ s'écrit $n = mq + r$ avec $q \geqslant 1$ et $r \in \{ 0, \ldots , m-1 \}$ ; ainsi, la sous-additivité permet d'écrire
$$
u(n) = u(mq + r) = u\big( m(q-1) + m+r \big) \leqslant (q-1) u(m) + u(m+r)
$$
puis
$$
\frac{u(n)}{n} \leqslant \frac{u(m)}{m} + \frac{1}{n} \max_{0\leqslant i<m} u(m+i)
$$
D'où, comme $u$ est contrôlée, en faisant tendre $n$, puis $m$ vers $\infty$, on obtient
$$
\limsup_{n\to\infty} \frac{u(n)}{n} \leqslant \liminf_{m\to\infty} \frac{u(m)}{m}
$$
autrement dit la suite $\big( u(n)/n \big)_{n\geqslant 1}$ converge. D'après le lemme \ref{lsa1}, sa limite est
$$
\phantom \qquad\qquad\qquad\qquad\qquad\qquad\qquad\qquad\quad \inf_{n \geqslant 1} \frac{u(n)}{n} \qquad\qquad\qquad\qquad\qquad\qquad\qquad\qquad\quad \square
$$

\subsection{Interversion infimum-supremum}

Il est question ici du fameux \emph{principle of the largest term} de \cite{LPS93}.

\begin{Le}\label{ls}
Si $\big(u_1(n)\big)_{n\geqslant 1}$, \ldots , $\big(u_r(n)\big)_{n\geqslant 1}$ sont $r$ suites à valeurs dans $[0, +\infty]$, on a l'égalité
\[
\limsup_{n\to\infty} \frac{1}{n} \log \sum_{i=1}^r u_i(n) = \max_{1\leqslant i\leqslant r} \limsup_{n\to\infty} \frac{1}{n} \log u_i(n)
\]
\end{Le}

\textbf{Démonstration :} De l'encadrement (on rappelle que les suites sont positives)
\[
\max_{1\leqslant i\leqslant r} u_i(n) \leqslant \sum_{i=1}^r u_i(n) \leqslant r \max_{1\leqslant i\leqslant r} u_i(n)
\]
on déduit que
\[
\limsup_{n\to\infty} \frac{1}{n} \log \sum_{i=1}^r u_i(n) = \limsup_{n\to\infty} \frac{1}{n} \log \max_{1\leqslant i\leqslant r} u_i(n)
\]
Puis
\begin{eqnarray*}
\limsup_{n\to\infty} \frac{1}{n} \log \max_{1\leqslant i\leqslant r} u_i(n) & = & \lim_{n\to\infty} \sup_{k\geqslant n} \frac{1}{k} \log \max_{1\leqslant i\leqslant r} u_i(k)\\
 & = & \lim_{n\to\infty} \max_{1\leqslant i\leqslant r} \left[\sup_{k\geqslant n} \frac{1}{k} \log u_i(k)\right]\\
 & = & \max_{1\leqslant i\leqslant r} \lim_{n\to\infty} \left[\sup_{k\geqslant n} \frac{1}{k} \log u_i(k)\right]\\
 & = & \max_{1\leqslant i\leqslant r} \limsup_{n\to\infty} \frac{1}{n} \log u_i(n)
\end{eqnarray*}
La troisième égalité vient du fait facile à vérifier que
\[
\max : [-\infty , +\infty]^r \rightarrow [-\infty , +\infty]
\]
est continue.\qed

\subsection{Transformation de Fenchel-Legendre}

Pour une présentation plus large, on renvoie à \cite{Mor67} et \cite{Roc70}. La preuve présentée ici est reprise de \cite[proposition 12.2.]{Cer07}.
Soit $(X, \tau)$ un espace vectoriel localement convexe. Notons $X^*$ son dual topologique. Définissons, pour $f : X \rightarrow [-\infty, +\infty]$, sa \emph{transformée de Fenchel-Legendre} par
$$
\forall \lambda \in X^* \qquad f^*(\lambda) := \sup_{x \in X}  \big( \langle \lambda | x \rangle - f(\lambda) \big)
$$
La transformée de Fenchel-Legendre de $g : X^* \rightarrow [-\infty, +\infty]$ se définit de façon analogue.

Etendons l'addition à $[-\infty , +\infty ]$ via
$$
\forall a \in \mathbb{R} \qquad \pm\infty \,\textup{\d{$+$}}\, a = \pm\infty \qquad \textrm{et} \qquad -\infty \,\textup{\d{$+$}}\, (+\infty) = -\infty
$$
ainsi que la multiplication par un réel via
$$
\begin{array}{ll} \forall \alpha > 0 & \alpha \cdot (\pm\infty) = \pm\infty \\
\forall \alpha < 0 & \alpha \cdot (\pm\infty) = \mp\infty
\end{array} \qquad \textrm{et} \qquad 0 \cdot (\pm\infty) = 0
$$
On dit que $f : X \rightarrow [-\infty, +\infty]$ est \emph{convexe} si
$$
\forall (x, y) \in X^2 \quad \forall \alpha \in [0, 1] \qquad f(\alpha x + (1-\alpha)y) \leqslant \alpha f(x) \,\textup{\d{$+$}}\, (1-\alpha) f(y)
$$
On notera que la seule fonction convexe prenant la valeur $-\infty$ est la fonction constante de valeur $-\infty$. Les fonctions convexes autres que les deux constantes $\pm \infty$ sont habituellement dénommées \emph{fonctions convexes propres}. On dit que $f : X \rightarrow [-\infty, +\infty]$ est \emph{concave} si $-f$ est convexe. On dit que $q : X \rightarrow [-\infty , +\infty ]$ est \emph{affine} si $q$ est convexe et concave. Les fonctions affines à valeurs dans $[-\infty, +\infty]$ sont alors les fonctions affines habituelles à valeurs dans $]-\infty, +\infty[$ et les deux fonctions constantes $\pm \infty$.

Le seul résultat qui nous intéresse ici est le suivant :
\begin{Pro}\label{fl}
Soit $(X, \tau)$ un espace vectoriel localement convexe et $f : X \rightarrow [-\infty , +\infty ]$. Alors
$$
f^{**} = f
$$
si et seulement si $f$ est convexe et semi-continue inférieurement relativement à la topologie faible $\sigma(X, X^*)$.
\end{Pro}

\textbf{Remarque :} Plus précisément, la transformation de Fenchel-Legendre réalise une bijection des fonctions convexes $\sigma(X, X^*)$-s.c.i. sur $X$ sur les fonctions convexes $\sigma(X^*, X)$-s.c.i. sur $X^*$, de sorte que, si $f$ est convexe et $\sigma(X, X^*)$-s.c.i., sa transformée de Fenchel-Legendre $f^*$ prend le nom de fonction \emph{convexe-conjuguée} de $f$.

\textbf{Démonstration :} L'implication directe découle du fait que, pour toute fonction $g : X^* \rightarrow [-\infty , +\infty]$, la fonction $g^*$ est convexe et $\sigma(X, X^*)$-s.c.i. Montrons la réciproque. On remarque que, si $\lambda \in X^*$, $\langle \lambda | \cdot \rangle - f^*(\lambda)$ est la plus grande fonction affine (y compris les deux fonctions affines $\pm \infty$) dirigée par $\lambda$ et inférieure à $f$. Il s'agit donc de voir que $f$ est la borne supérieure de l'ensemble des fonctions affines continue plus petites que $f$ (y compris les deux fonctions affines $\pm \infty$). Si $f = \pm \infty$, le résultat est immédiat. Sinon, $f$ est une fonction convexe propre. Définissons l'épigraphe de $f$ par $epi(f) = \{ (x, t) \in X \times \mathbb{R} ; f(x) \leqslant t \}$. Le fait que $f$ soit convexe et $\sigma(X, X^*)$-s.c.i. assure que $epi(f)$ est convexe et fermé relativement à $\tau$. Soit $(x, t) \in X \times \mathbb{R} \setminus epi(f)$. Le théorème de Hahn-Banach dans l'espace localement convexe $X \times \mathbb{R}$ donne l'existence d'un hyperplan fermé séparant strictement $(x, t)$ et $epi(f)$. Si cet hyperplan n'est pas vertical, il correspond à une fonction affine plus petite que $f$. Sinon, l'hyperplan est de la forme $H \times \mathbb{R}$ où $H$ est un hyperplan affine fermé de $X$ et $f(x) = +\infty$. Soit alors $p$ une fonction affine continue sur $X$ telle que $p(x) > 0$ et $p(y) = 0$ pour tout $y \in H$. Comme $f$ est une fonction convexe propre, il existe une fonction affine continue et finie $q$ inférieure à $f$ : en effet, il existe $x \in X$ tel que $f(x) \in ]-\infty , +\infty[$ et un hyperplan fermé $H$ séparant $(x, f(x) - 1)$ de $epi(f)$ ; cet hyperplan $H$ n'est pas vertical et correspond à une fonction affine continue et finie $q$ inférieure à $f$. Ainsi, pour tout $\alpha > 0$, $q + \alpha p$ est toujours une fonction affine continue inférieure à $f$. Il suffit alors de choisir $\alpha$ tel que $q(x)+\alpha p(x) > t$.\qed

\newpage


\section{Compléments}

\subsection{PGD fort}

Il existe des conditions analogues à celles du cas i.i.d. donnant le PGD fort à partir du PGD faible. La plus simple est le cas où il existe un ensemble convexe mesurable $K$ d'adhérence compacte tel que $(\mu_\Lambda)_{\Lambda \in \mathbb{B}}$ soit porté par $K^{\mathbb{Z}^d}$. Voici une condition plus faible (et classique dans le cas i.i.d.) :
\begin{Th}[PGD fort]
Soit $(\mu_\Lambda)_{\Lambda \in \mathbb{B}}$ un s.c.i.t. a.d.i. et c.l. On suppose qu'il existe $K$ convexe mesurable relativement compact de $X$ tel que, en notant $M_K$ sa jauge,
$$
\limsup_{n \to \infty} \frac{1}{|\Lambda(n)|} \log \mathbb{E}\left( e^{|\Lambda(n)| M_K \big( \mathfrak{m}_{\Lambda(n)}(\sigma) \big) } \right) < + \infty
$$
alors $(\mu_n)_{n \geqslant 1}$ vérifie un PGD fort.
\end{Th}

\textbf{Démonstration :} Soit $M < 0$ et $t = \overline{p}(M_K) - M$. L'inégalité de Tchebychev permet d'obtenir
\begin{align*}
\limsup_{n \to \infty} \frac{1}{|\Lambda(n)|} \log \mu_n \big( X \setminus (tK) \big)
 &\leqslant \limsup_{n \to \infty} \frac{1}{|\Lambda(n)|} \log \mu_n \big( |\Lambda(n)| M_K > |\Lambda(n)| t \big) \\
 &\leqslant \limsup_{n \to \infty} \frac{1}{|\Lambda(n)|} \log \left( e^{-|\Lambda(n)| t} \int e^{|\Lambda(n)| M_K} d \mu_n \right) = \overline{p}(M_K) - t = M\hspace{-5pt}
\end{align*}
On en déduit que la suite $(\mu_n)_{n \geqslant 1}$ est exponentiellement tendue, ce qui permet de conclure. Voici tout de même la preuve. La borne inférieure pour $X \setminus K$ donne
$$
\sup_{X \setminus \overline{K}} s \leqslant t
$$
donc $\{ x \in X \, ; \, s(x) > t \} \subset \overline{K}$ et $s$ est coercive. Soit maintenant $A \in \mathcal{F}$. Soit $M < 0$. Il existe $K \in \mathcal{F}$ relativement compact tel que $\overline{p}(\delta_{X \setminus K}) \leqslant M$. Alors
$$
\mu_n(A) \leqslant \mu_n(A \cap K) + \mu_n(X \setminus K)
$$
La borne supérieure faible pour $A \cap K$ et le lemme \ref{ls} montrent que
$$
\limsup_{n \to \infty} \frac{1}{|\Lambda(n)|} \log \mu_n(A) \leqslant \max \big( M, \sup_{\overline{A \cap K}} s \big) \leqslant \max \big( M, \sup_{\overline{A}} s \big)
$$
On obtient la borne supérieure en faisant tendre $M$ vers $-\infty$.\qed

\subsection{Cas d'égalité des entropies inférieure et supérieure}

En récrivant le corollaire \ref{lsacor}, on a montré que, pour tout $x \in X$, pour tout $C \in \mathcal{C}_x$ et pour tout $\varepsilon \in ]0, 1[$,
$$
\liminf_{n \to \infty} \frac{1}{|\Lambda(n)|} \log \mu_n(x+C) \geqslant \limsup_{m \to \infty} \frac{1}{|\Lambda(m)|} \log \mu_m \big( x + (1-\varepsilon)C \big)
$$
On donne ici une condition suffisante pour que
$$
\lim_{\varepsilon \to 0^+} \limsup_{m \to \infty} \frac{1}{|\Lambda(m)|} \log \mu_m \big( x + (1-\varepsilon)C \big) = \limsup_{m \to \infty} \frac{1}{|\Lambda(m)|} \log \mu_m(x + C)
$$

\begin{Pro}
On suppose qu'il existe $y \in C$, $\beta \in [0, 1[$ et $\alpha_0 > 0$ tels que, pour tout sous-ensemble fini $S$ de $\mathbb{Z}^d \setminus \{ 0 \}$ et pour tout $D \in \mathcal{C}(X^S)$,
$$
\mathbb{P} \big( \sigma(0) \in y + \beta(C - y) ; \eta_S \in D \big) \geqslant \alpha_0 \mathbb{P} \big( \sigma_S \in D \big)
$$
Alors
$$
\liminf_{n \to \infty} \frac{1}{|\Lambda(n)|} \log \mu_n(x+C) = \limsup_{m \to \infty} \frac{1}{|\Lambda(m)|} \log \mu_m(x + C)
$$
\end{Pro}

\textbf{Remarque :} La condition qui apparaît ici est plus forte que l'hypothèse de contrôle local pour $V$. Elle impose, essentiellement, que la loi de $\sigma(0)$ donne de la masse à l'ouvert convexe $C$, conditionnellement au reste de la configuration.

\textbf{Démonstration :} On définit des jauges tronquées de $V$ : pour tout $\beta \in [0, 1[$,
\[
M_{V, \beta} = \frac{(\beta \vee M_V) - \beta}{1-\beta}
\]
et on note $\phi_\beta$ la fonction convexe $M_{V, \beta}(\cdot - y)$.

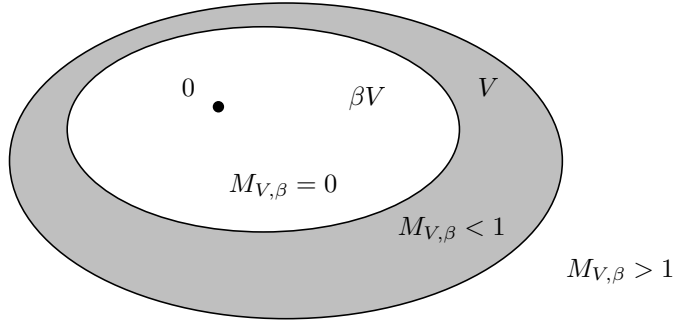
\begin{center}
\def\JPicScale{0.7}
\input{jaugetronquee.pst}
\end{center}

Par un cheminement analogue à celui de la démonstration du lemme sous-additif \ref{lsa}, on obtient
\begin{align*}
\frac{1}{|\Lambda(n)|}& \log \mu_m \big( x + (1-\varepsilon)C \big)\\
 &=\frac{1}{|\Lambda(n)|} \log \mathbb{P}\left(\phi_\beta\big(\mathfrak{m}_{\Lambda(n)} (\sigma)\big) < 1 -\frac{\varepsilon}{1-\beta} \right)\\
 &\geqslant \frac{1}{|\Lambda(m)|} \log \mathbb{P}\big(\phi_\beta\big(\mathfrak{m}_{\Lambda(m)} (\sigma)\big) < 1 - \frac{2\varepsilon}{1-\beta} \big) - \frac{c(m)}{|\Lambda(m)|} + \rho_{m, n} \log \alpha_0\\
\end{align*}
En faisant tendre $n$ vers $\infty$, puis $\varepsilon$ vers $0^+$, et enfin $m$ vers $\infty$, on obtient
$$
\lim_{\varepsilon \to 0^+} \limsup_{m \to \infty} \frac{1}{|\Lambda(m)|} \log \mu_m \big( x + (1-\varepsilon)C \big) = \limsup_{m \to \infty} \frac{1}{|\Lambda(m)|} \log \mu_m(x + C)
$$
ce qu'on voulait.\qed

\medskip

La démonstration s'adapte aux limites projectives en considérant l'indice $i \in J$ tel que $C = f_i^{-1}(C_i)$.

\subsection{Propriétés de la pression}

Etendons l'addition à $[-\infty , +\infty ]$ via
$$
\forall a \in \mathbb{R} \qquad \pm\infty \,\textup{\d{$+$}}\, a = \pm\infty \qquad \textrm{et} \qquad -\infty \,\textup{\d{$+$}}\, (+\infty) = -\infty
$$
ainsi que la multiplication par un réel via
$$
\begin{array}{ll} \forall \alpha > 0 & \alpha \cdot (\pm\infty) = \pm\infty \\
\forall \alpha < 0 & \alpha \cdot (\pm\infty) = \mp\infty
\end{array} \qquad \textrm{et} \qquad 0 \cdot (\pm\infty) = 0
$$
On dit que $f : X \rightarrow [-\infty, +\infty]$ est \emph{convexe} si
$$
\forall (x, y) \in X^2 \quad \forall \alpha \in [0, 1] \qquad f(\alpha x + (1-\alpha)y) \leqslant \alpha f(x) \,\textup{\d{$+$}}\, (1-\alpha) f(y)
$$

\begin{Pro}
La pression est une fonction convexe.
\end{Pro}

\textbf{Démonstration :} Pour tout $n \geqslant 1$, l'inégalité de Hölder permet de montrer que
$$
\lambda \in X^{*m} \mapsto \frac{1}{v_n} \log \int e^{v_n \lambda} d\mu_n \in [-\infty , +\infty]
$$
est convexe. On conclut en remarquant qu'une limite supérieure de fonctions convexes est convexe.\qed

\medskip

Comme dans le cas indépendant, on montre :
\begin{Pro}
Si $\mu_{\Lambda(1)}$ est convexe-tendue localement, alors la pression est semi-continue inférieurement relativement à $\sigma(X, X^*)$.
\end{Pro}

\subsection{Champs contrôlés localement}

Le modèle d'Ising (et plus généralement les mesures de Gibbs à espace d'états borné) est a.d.i. avec $g(1) = 0$. On montre par exemple que, pour toute boîte $\Lambda \subset \mathbb{Z}^d \setminus \{ 0 \}$ et pour tout $A \in \mathcal{F}_\Lambda$,
$$
\mathbb{P}(\sigma(0) = 1 ; \sigma_\Lambda \in A) \geqslant \gamma \mathbb{P}(\sigma(0) = 1) \mathbb{P}(\sigma_\Lambda \in A)
$$
avec
$$
\gamma = \frac{1}{(1+e^{8\beta})\mathbb{P}(\sigma(0) = 1)} > 0
$$


\bibliographystyle{alpha-fr}
\bibliography{../cramer}

\end{document}